\documentclass{amsart}
\usepackage{amscd, amssymb, portland, calc, ifthen}
\usepackage[matrix,arrow,curve,cmtip]{xy}
\CompileMatrices

\newtheoremstyle{mythm}{}{}%
  {\itshape}%	Body font
  {}%         	Indent amount (empty = no indent, \parindent = para indent)
  {\bfseries}% 	Thm head font
  {}%        	Punctuation after thm head
  { }%     	Space after thm head (\newline = linebreak)
  {\thmnumber{#2.\hspace{1.5mm}}\thmname{#1}\thmnote{ #3}.}%  Thm head spec

\newtheoremstyle{myrmk}{}{}%
  {}%		Body font
  {}%         	Indent amount (empty = no indent, \parindent = para indent)
  {\bfseries}% 	Thm head font
  {}%        	Punctuation after thm head
  { }%     	Space after thm head (\newline = linebreak)
  {\thmnumber{#2.\hspace{1.5mm}}\thmname{#1}\thmnote{ #3}.}%  Thm head spec

\numberwithin{equation}{subsection}

\theoremstyle{mythm}
\newtheorem{thm}[subsection]{Theorem}
\newtheorem{prop}[subsection]{Proposition}
\newtheorem{lemma}[subsection]{Lemma}
\newtheorem{cor}[subsection]{Corollary}

\theoremstyle{myrmk}

\newtheorem{const}[subsection]{Construction}

\theoremstyle{plain}
\newtheorem*{thm*}{Theorem}
\newtheorem*{prop*}{Proposition}
\newtheorem*{cor*}{Corollary}

\theoremstyle{definition}
\newcommand{\mathbold}{\bf}

\newcommand{\sep}{^{\mathrm{sep}}}

\def\longisomap{{\buildrel \sim\over\longrightarrow}} % an isomorphism

\newcommand{\subsec}[1]{{\begin{trivlist}\item\em\large#1\end{trivlist}}}

\newcommand{\longlabelmap}[1]{{\,\buildrel #1\over\longrightarrow\,}}

\newcommand{\longmap}{{\,\longrightarrow\,}}

\newcommand{\Hom}{{\mathrm{Hom}}}

\newcommand{\res}{{\mathrm{res}}}
\DeclareMathOperator{\dlog}{dlog}

\newcommand{\Spec}{{\mathrm{Spec}}}
\newcommand{\Mod}{{\:\mathrm{mod}\:}}
\newcommand{\coker}[1]{{{\mathrm{coker}}\left[#1\right]}}
\newcommand{\kernel}[1]{{{\ker}\left[#1\right]}}

\newcommand{\Gal}{{\mathrm{Gal}}}

\newcommand{\gr}{{\mathrm{gr}}}
\newcommand{\grs}{{\mathrm{gr}}\,}

\newcommand{\ff}{{\mathbold F}}
\newcommand{\qq}{{\mathbold Q}}

\newcommand{\zz}{{\mathbold Z}}

\newcommand{\nn}{{\mathbold N}}
\newcommand{\G}{{\mathbold G}}
\newcommand{\et}{{\mbox{\scriptsize{\'et}}}}

\newcommand{\textfildot}{{\scriptscriptstyle\bullet\!}}
\newcommand{\scriptdot}{{^{_{_{_{\textfildot}}}}}}

\newcommand{\m}{{\mathfrak{m}}}
\newcommand{\card}[1]{\,|#1|\,}
\newcommand{\codim}{{\mathrm{codim}\,}}

\newcommand{\CW}{{{CW}}}
\def\longisomap{{\,\buildrel \sim\over\longrightarrow\,}} % an isomorphism
\newcommand{\proot}[1]{{#1^{p^{-1}}}}
\newcommand{\comment}[1]{}

\newcommand{\filreg}[1]{{\fil'_{#1}}}
\newcommand{\fillog}{{\fil}}

\newcommand{\gc}{^{\mathrm{g}}}

\newcommand{\KM}[1]{{K^{\mathrm{M}}_{#1}}}
\newcommand{\KMb}[1]{{k_{#1}}}
\newcommand{\fil}{{\mathrm{fil}}}
\newcommand{\del}{{\partial}}
\newcommand{\arnv}[1]{{\mathrm{{ar_{\scriptstyle\mathrm{n}}^{#1}}}}}

\newcommand{\ark}{{\mathrm{ar_{\scriptscriptstyle K}}}}
\newcommand{\swk}{{\mathrm{sw_{\scriptscriptstyle K}}}}
\newcommand{\art}{{\mathrm{ar}}}
\newcommand{\AS}{{\xi}}

\newcommand{\cech}[1]{{\check{#1}}}
\newcommand{\eh}{{\hat{e}}}

\newcommand{\Omlog}[2]{{\omega^{{#1}}_{{{#2}}}}}
\newcommand{\p}{\mathfrak{p}}

\renewcommand{\m}{\p}
\renewcommand{\leq}{\leqslant}
\renewcommand{\geq}{\geqslant}

\makeatletter
\def\@seccntformat#1{\@ifundefined{#1@cntformat}%
{\csname the#1\endcsname\quad}%  default
{\csname #1@cntformat\endcsname}% individual control
}
\def\section@cntformat{\thesection.\enspace}
\def\subsection@cntformat{\thesubsection.}
\makeatother

\begin{document}

\newcommand{\kk}{{k}}
\newcommand{\kkt}{{\tilde{k}}}
\newcommand{\kkgc}{{k\gc}}

% Balint Virag's macros
%%%%%%%%%%%%%%%%%%%%%%%
%\newcommand\mnote[1]{\marginpar{\ \\ \small \tt #1}}
\newcommand\mnote[1]{}
\newcommand\bel[1]{{\mnote{#1}}\begin{equation}\label{#1}}
\newcommand\lb[1]{\label{#1}\mnote{#1}}

% Time
\newcounter{hour}\newcounter{minute}
\newcommand{\printtime}{\setcounter{hour}{\time/60}%
	\setcounter{minute}{\time-\value{hour}*60}%
	\ifthenelse{\value{hour}<10}{0}{}\thehour:%
	\ifthenelse{\value{minute}<10}{0}{}\theminute}

\setlength{\marginparwidth}{3cm}

\title{Kato's conductor and generic residual perfection}
\author[J.~Borger]{James M. Borger}
\address{University of Chicago, Chicago, Illinois, USA}
%\curraddr{}
\email{borger@math.uchicago.edu}
%\thanks{To A,T,E, and I.}
%\date{\today. \printtime}
\begin{abstract}
  Let~$A$ be a complete discrete valuation ring with possibly
  imperfect residue field, and let~$\chi$ be a rank-one Galois
  representation over~$A$.  I show that the non-logarithmic variant of
  Kato's Swan conductor is the same for~$\chi$ and the pullback
  of~$\chi$ to the generic residual perfection of~$A$.  This implies
  the conductor from ``Conductors and the moduli of residual
  perfection''~\cite{Borger:Cond-moduli} extends the non-logarithmic
  variant of Kato's.
\end{abstract}
\maketitle

\section*{Introduction}
Let~$A$ be a complete (or even henselian)
discrete valuation ring with fraction field~$K$.
When the residue field~$\kk$ of~$A$ is perfect, there is~\cite[IV,VI]{Serre:CL}
a well-known and satisfactory theory of wild ramification
over~$A$.  We understand much less, though, when we do not require
that~$\kk$ be perfect.  In this context, Kato~\cite{Kato:Imp} has
developed a good abelian theory: to every one-dimensional Galois
representation~$\chi$ over~$A$, he gives a non-negative integer, the
Kato-Swan conductor~$\swk(\chi)$ of~$\chi$, that measures the extent
to which~$\chi$ is wildly ramified.  It is the
pole order, in the logarithmic sense, of a certain differential
form---his refined Swan conductor.  We can just as well, however,
consider the order in the usual, non-logarithmic sense.  I will call
this the Kato-Artin conductor of~$\chi$ and denote it by~$\ark(\chi)$.
When~$\kk$ is perfect, it agrees with the usual Artin conductor,
and~$\swk(\chi)$ agrees with the usual Swan conductor.

There have been several recent
attempts~\cite{Abbes-Saito:Ram-I,Boltje.et.al:cond,Borger:Cond-moduli,Zhukov:Ram}
to give a general approach to non-abelian wild ramification over such
rings~$A$.  In one of them, I showed how to associate a non-negative
integer~$\art(\rho)$ to any Galois representation~$\rho$ over~$A$.
The purpose of this paper is to show that when~$\rho$ is
one-dimensional,~$\art(\rho)$ agrees with~$\ark(\rho)$.

Let~$A\gc$ be the generic residual perfection~\cite{Borger:Cond-moduli} of~$A$.
It is a residually perfect complete discrete valuation ring 
of ramification index one over~$A$ and is, in a
certain sense, universally generic with respect to these properties.
Let~$K\gc$ be its fraction field.  Most of this paper is devoted to the
proof of the following result.

\begin{thm*}
If~$\chi$ is a class in~$H^1(K,\qq/\zz)$ and~$\chi'$ is its
image in~$H^1(K\gc,\qq/\zz)$, then~$\ark(\chi)=\ark(\chi')$,
%Let~$\chi$ be a class in~$H^1(K,\qq/\zz)$, and let~$\chi'$ be its
%image in~$H^1(K\gc,\qq/\zz)$.  Then~$\ark(\chi)=\ark(\chi')$,
\end{thm*}

The intuitive reason why this should be true is that
% the Kato-Artin conductor measures the pole order of a
% certain differential form and
the order of any differential form on~$A$ should remain
unchanged when the form is pulled back to~$A\gc$.  When~$A$ is of equal
characteristic, this has meaning and, once the necessary foundations
are laid, is essentially a proof.  In fact, the observation that there
are residually perfect extensions with this property is what led to
the definition~\cite{Borger:Cond-moduli} of the general Artin
conductor.  In mixed characteristic, however, this
provides little more than motivation, and most of this paper
is spent pushing it through to a real proof.

In section~\ref{sec-kato-theory}, I recall Kato's theory, prove some
basic results, and define the Kato-Artin conductor.  The proof of the
theorem when~$A$ is of equal characteristic is in
section~\ref{sec-equal-kato}.  It uses Matsuda's
refinement~\cite{Matsuda:Swan} of Kato's refined Swan conductor.
Because the proof in equal characteristic is so much simpler than the
proof in mixed characteristic, I encourage the reader to read it
first.  The basic technique in mixed characteristic is to use Kato's
description~\cite[4.1]{Kato:Imp} (following
Bloch-Kato~\cite{Bloch-Kato:p-adic}) of certain graded pieces of
cohomology groups in terms of explicit $K$-theoretic symbols and then
understand how these symbols behave under pullback to~$A\gc$.
Section~\ref{sec-lemmas} contains a commutative diagram that encodes
this behavior, and section~\ref{sec-mixed-kato} gives the proof in
mixed characteristic.  In the final section, I show how the theorem
implies that~$\art(\rho)$ and~$\ark(\rho)$ agree for one-dimensional
Galois representations~$\rho$.

\section*{Conventions}
If~$A$ is a discrete valuation ring,~$\p_A$ will denote its maximal
ideal, and~$U^{\scriptdot}_A$ will denote the filtration of~$A^*$
with~$U^0_A=A^*$ and~$U^i_A=1+\m_A^i$ for positive integers~$i$.  If
its residue field~$A/\m_A$ is perfect and has positive characteristic,
then~$s_A$ will denote the unique multiplicative section~\cite[II \S 4
Prop.~8]{Serre:CL} of the reduction map~$A\to A/\p_A$.  If~$x$ is an
element of~$A$, its image in~$A/\p_A$ will be denoted~$\bar{x}$.  An
extension of~$A$ is a discrete valuation ring~$B$ equipped with an
injective local homomorphism~$A\to B$.  We will denote its
ramification index by~$e_{B/A}$.

Throughout,~$A$ will denote a henselian~\cite{Raynaud:Hens}
discrete valuation ring, held
fixed within each subsection, with fraction field~$K$ and residue
field~$k$. We will assume for simplicity of exposition that~$k$ always has
characteristic~$p$, where~$p$ is a fixed prime number.  We will
say~$A$ is {\em of equal characteristic} if~$K$ has characteristic~$p$
and is {\em of mixed characteristic} if~$K$ has characteristic~$0$.

\section{The generic residual perfection}

The purpose of this section is to recall the main
facts~\cite[1.14, 2.4]{Borger:Cond-moduli} about the generic
residual perfection we will need later.
%\marginpar{ref} %% Unstable reference

\subsection{}
The generic residual perfection of (the completion of)~$A$
is the~$A$-algebra~$A\gc$
corresponding to the generic point of the moduli space of residual
perfections of~$A$.  It is a complete discrete
valuation ring with~$e_{A\gc/A}=1$ whose residue field~$k\gc$
is perfect. %~\cite[1.15]{Borger:Cond-moduli}.
Its fraction field will be denoted~$K\gc$.

\begin{thm} %[\cite{Borger:Cond-moduli}]
\lb{thm-expl}
Let~$T$ be a lift to~$A$ of a $p$-basis of its residue field~$\kk$,
and let~$\pi\in A$ be a uniformizer.  For each element~$t\in T$,
let~$u_{t,1}, u_{t,2},\dots\in\kkgc$ be the unique sequence such that
the image of~$t$ in~$A\gc$
is~$s_{A\gc}(\bar{t})+s_{A\gc}(u_{t,1})\pi+s_{A\gc}(u_{t,2})\pi^2+\cdots$.
Let~$R$ be the free polynomial algebra~$\kk[T\times \zz_{>0}]$.  Then
the map~$R\to\kkgc$ determined by~$(t,j)\mapsto u_{t,j}$ induces an
isomorphism from the fraction field of~$R^{p^{-\infty}}$ to~$\kkgc$.
\end{thm}

\begin{thm} %[\cite{Borger:Cond-moduli}]
\lb{thm-Galois-surj}
Fix a separable closure of~$K\gc$.  Then the map~$G_{K\gc}\to G_K$ of
the corresponding absolute Galois groups is surjective.  The induced
maps of inertia groups and wild inertia groups are also surjective.
\end{thm}

\section{Kato's theory}
\lb{sec-kato-theory}

The purpose of this section is to collect some results in
Kato's theory~\cite{Kato:Imp}.  Let us first recall the basics.

\subsection{}
\lb{sbsc-Kato-notation}
Let~$F$ be a field and let~$n>0$ and~$r$ be integers.  If~$n$ is
invertible in~$F^*$,
let~$\zz/n(r)$ be the $r$-th Tate twist of the constant sheaf~$\zz/n$
on the \'etale topology (of Grothendieck~\cite{Grothendieck:Etale-top})
of~$F$. If the characteristic
of~$F$ is~$p>0$, write~$n=mp^s$, where~$p\nmid m$, and let~$\zz/n(r)$
be the complex
\[
\zz/m(r)\oplus W_s\Omega^r_{F,\log}[-r]
\]
of abelian sheaves on~$\Spec(F)_\et$.  Here,~$W_s\Omega^r_{F,\log}$ is
the piece of degree~$r$ of the logarithmic part \cite[I
5.7]{Illusie:dRW} of Deligne and Illusie's de\! Rham-Witt
complex~$W_s\Omega_{F}^\scriptdot$ on~$\Spec(F)_{\et}$.

For positive integers~$q$, write~$H^q_n(F)=H^q(F,\zz/n(q-1))$, and
let~$H^q(F)$ be the colimit of~$H^q_n(F)$ over the integers~$n$
(ordered by divisibility).  The natural map~$H^q_n(F)\to H^q(F)$ is an
isomorphism of~$H^q_n(F)$ with the $n$-torsion
of~$H^q(F)$. I will usually identify the two without comment.

%\subsection{}
Let
\[
h_F:F^*\longmap H^1(F,\zz/n(1))
\]
be the connecting homomorphism of the
Kummer triangle
\[
\zz/n(1)\longmap\G_m\longlabelmap{n}\G_m\longmap\zz/n(1)[1].
\]
(When~$n$ is a power of the characteristic of~$F$, the existence of
such a triangle follows from the theory of the de\! Rham-Witt
complex~\cite[I 3.23.2, I 5.7.1]{Illusie:dRW}.) 
Let~$\KM{r}(F)$ denote the~$r$-th Milnor
$K$-group~\cite{Milnor:K-theory} of~$F$.  Then there is a homomorphism
\[
K^\mathrm{M}_r(F) \longmap H^{r}(F,\zz/n(r)),
\]
also denoted~$h_F$, sending~$\{x_1,\dots,x_r\}$ (for~$x_1,\dots,x_r\in
F^*$) to the cup product
\[
h_F(x_1)\cup\cdots\cup h_F(x_r).
\]
For~$\chi\in H^q_n(F)$, let~$\{\chi,x_1,\dots,x_r\}$ denote~$\chi\cup
h_F(\{x_1,\dots,x_r\})$.  Taking the colimit over integers~$n$, we
get a pairing
\[
H^q(F)\otimes K^\mathrm{M}_r(F)\longmap H^{q+r}(F).
\]
We will use~$\{\chi,x_1,\dots,x_r\}$ to denote the image
of~$\chi\otimes\{x_1,\dots,x_r\}$ under this pairing as well.

If the characteristic of~$F$ is~$p$,
let~$\AS_s:W_s\Omega^{q-1}_{F}\twoheadrightarrow H^q_{p^s}(F)$
denote the higher Artin-Schreier maps~\cite[1.3]{Kato:Imp}.

\subsection{}
For any non-negative integer~$n$, let $\fillog_n H^q(K)$ be the
subgroup of classes~$\chi$ that have the
property~$\{\chi|_{A'},1+\m_A^n\m_{A'}\}=0$ for every henselian
extension~$A'$ of~$A$. This gives~\cite[2.2, 6.3]{Kato:Imp} an
exhaustive increasing filtration of~$H^q(K)$.
The {\em Kato-Swan conductor}
(or {\em logarithmic Kato conductor})  of a class~$\chi\in H^q(K)$
is the smallest integer~$n$ such that~$\chi\in\fillog_n H^q(K)$.  It
is denoted~$\swk(\chi)$.

\subsection{}
\lb{sbsc-tilde}
Let~$\tilde{A}$ denote the henselization of the localization of the
polynomial algebra~$A[T]$ at the ideal generated by~$\m_A$.
Then~$\tilde{A}$ is a henselian discrete valuation ring, and
for any uniformizer~$\pi$ of~$A$, we have~\cite[6.3]{Kato:Imp}
\[
\chi\in\fillog_n H^q(K) \ \Longleftrightarrow
\ \{\chi|_{\tilde{A}},1+\pi^{n+1}T\}=0.
%\fillog_n H^q(K) = \left\{\chi\in H^q(K)\mid\{\chi|_{\tilde{A}},1+\pi^{n+1}T\}=0\right\}.
\]
We will denote the fraction field of~$\tilde{A}$ by~$\tilde{K}$
and the residue field by~$\tilde{k}$.

\subsection{}
\lb{sbsc-fil0-str}
The map~$H^q(k)\to H^q(K)$ extends naturally to an exact sequence
\[
0 \longmap H^q(k) \longmap \fil_0 H^q(K) \longmap H^{q-1}(k) \longmap 0.
\]
Given a uniformizer~$\pi$ of~$A$, the map~$\psi \mapsto \{\psi,\pi\}$
is a splitting~\cite[6.1]{Kato:Imp}.

\subsection{}
\lb{sbsc-log-diff}
The reduction map
\[
(A-\{0\})/U^1_A \longmap \kk
\]
gives~$\kk$ the structure of a log ring~\cite[1.1]{Kato:Log}.
Let~$\Omlog{1}{\kk}=\Omlog{1}{\kk/\zz}$ denote
the~$\kk$-module of absolute K\"ahler differentials with respect
to this log structure~\cite[1.7]{Kato:Log}.  For~$q\in\nn$,
let~$\Omlog{q}{\kk}$ denote the $q$-th exterior power
of~$\Omlog{q}{\kk}$.
There is a natural exact sequence
\[
0 \longmap \Omega^q_{\kk} \longmap \Omlog{q}{\kk}
  \longlabelmap{\res} \Omega^{q-1}_{\kk} \longmap 0.
\]
(Of course,~$\Omega^{\scriptdot}_{\kk}$
means~$\Omega^{\scriptdot}_{\kk/\zz}$.) Given a uniformizer~$\pi$
of~$A$, the map~$\eta\mapsto\eta\wedge\dlog(\pi)$ is a splitting.

\subsection{}
There is a unique map~$\lambda_A:\Omlog{q-1}{\kk}\longmap \fil_0
H^q_p(K)$ that gives rise to a map of sequences
\[\xymatrix{
0\ar[r] & \Omega^{q-1}_{\kk}\ar[r]\ar@{>>}[d]_{\AS_1} & {\Omlog{q-1}{\kk}}\ar[r]\ar@{>>}[d]_{\lambda_A} & \Omega^{q-2}_{\kk}\ar[r]\ar@{>>}[d]_{\AS_1} & 0 \\
0\ar[r] & H^q_p(\kk)\ar[r] & {\fil}_0 H^q_p(K)\ar[r] & H^{q-1}_p(\kk)\ar[r] & 0
}
\]
respecting the splittings in~\ref{sbsc-fil0-str}
and~\ref{sbsc-log-diff} (for every uniformizer~$\pi$ of~$A$).

The following theorem \cite[5.1, 5.2, 5.3]{Kato:Imp} is fundamental.

\begin{thm}
\lb{thm-kato-thy}
%Let~$A$ be a henselian discrete valuation ring, and 
Let~$n$ be a positive integer.  Then for any class~$\chi\in\fil_n
H^q(K)$, there is a unique
element~$\eta\in\m_A^{-n}\otimes_A\Omlog{q}{\kk}$ such that for every
henselian extension~$A'$ of~$A$ and every element~$z$ of~$\m_A^n A'$, we have
\[
\{\chi,1+z\}=\lambda_{A'}(z\eta).
\]
Furthermore, the function~$\chi\mapsto\eta$ induces an injective homomorphism
\[
\kappa_n:\gr_n H^q(K)\longmap\m_A^{-n}\otimes\Omlog{q}{\kk}.
\]
\end{thm}
I will typically write~$\kappa_n(\chi)$ for the image under~$\kappa_n$
of the graded class of~$\chi$.  It is
called the {\em refined Swan conductor of~$\chi$}.

\subsection{}
Let us now consider the non-logarithmic analogue of the Kato-Swan
conductor. Let~$\chi$ be a class in~$H^1(K,\qq/\zz)$ and
put~$n=\swk(\chi)$.  Define
\[
\ark(\chi) = \left\{
\begin{array}{ll}
0   & \mbox{if~$\chi$ is unramified} \\
1   & \mbox{if~$\chi$ is tame and ramified} \\
n   & \mbox{if~$\chi$ is ramified
        and~$\kappa_n(\chi)\in\m_A^{-n}\otimes{\Omega}^1_{\kk}$} \\
n+1 & \mbox{if~$\chi$ is ramified
        and~$\kappa_n(\chi)\not\in\m_A^{-n}\otimes{\Omega}^1_{\kk}$}
\end{array}\right.
\]
We call~$\ark(\chi)$ the {\em Kato-Artin conductor} of~$\chi$.  It is
the natural non-logarithmic analogue of Kato's Swan conductor.  As
mentioned in the introduction, when~$\swk(\chi)$ is not
zero,~$\ark(\chi)$ can be viewed as the order of the pole
of~$\kappa_n(\chi)$ in
the usual sense and~$\swk(\chi)$ can be viewed as the order in the
logarithmic sense.

Matsuda~\cite[3.2.5]{Matsuda:Swan} has used what is essentially the
same conductor.  His is one less than~$\ark(\chi)$ except when~$\chi$
is unramified, in which case both are zero.

\subsec{Basic facts}

The rest of this section contains some propositions we will need in
the proof of the theorem in mixed characteristic.  All the proofs
are straightforward.

\begin{prop}
\lb{prop-p-div-sw}
Let~$\chi$ be a class in~$H^1(K)$.  If~$\ark(\chi)$ and~$\swk(\chi)$
have the same value, then it is a multiple of~$p$.
\end{prop}
\begin{proof}
\cite[5.4]{Kato:Imp}
\end{proof}

Let~$A'$ be a finite extension
of~$A$ of ramification index~$e$, let~$K'$ denote its fraction field,
and let~$k'$ denote its residue field.  Let~$n$,~$q$, and~$s$ be
positive integers.

\begin{prop}
\lb{prop-rsw-natural}
%For all positive integers~$q$ and~$n$,
The following diagram commutes:
\[\xymatrix{
{\gr}_{en} H^q(K')\ar[r]^{\kappa_{en}}
   & {\m}_{A'}^{en}\otimes_{A'}\Omlog{q}{k'} \\
{\gr}_{n} H^q(K)\ar[r]^{\kappa_n}\ar[u]
   & {\m}_{A}^{n}\otimes_A\Omlog{q}{\kk}\ar[u].
}\]
\end{prop}
\begin{proof}
Use the uniqueness statement in~\ref{thm-kato-thy}.
\end{proof}

\begin{cor}
\lb{cor-tame-bc}
If the extension~$A'/A$ is tame, then
for 
%any integer~$q\geq 1$ and
any class~$\chi\in H^q(K)$,
we have~$\swk(\chi|_{A'}) = e\,\swk(\chi)$
\end{cor}
\begin{proof}
The map~$\Omega^{\scriptdot}_k \to \Omega^{\scriptdot}_{k'}$ is injective.
\end{proof}

\begin{prop}
\lb{prop-G-inv}
If the extension~$K'/K$ is Galois with group~$G$ and its degree is not
a multiple of~$p$, the natural maps
\begin{align*}
H^q_{p^s}(K) & \longmap  H^q_{p^s}(K')^G \\
\fil_{n-1} H^q_{p^s}(K) & \longmap  \fil_{en-1}H^q_{p^s}(K')^G
\end{align*}
are isomorphisms.% for every positive integer~$n$.
%all integers~$n\geq 1,q\geq 1,s\geq 0$.
\end{prop}
\begin{proof}
Since the order of~$G$ is relatively prime to~$p$, the
groups~$H^i(G,H^j_{p^s}(K'))$ are zero for~$i>0$.
The existence of a spectral sequence
\[
H^i(G,H^j_{p^s}(K')) \Rightarrow H^{i+j}_{p^s}(K),
\]
implies the first map is an isomorphism.  The second then 
is by~\ref{cor-tame-bc}.
\end{proof}

\begin{prop}
\lb{prop-fil0-structure}
%For all positive integer~~$q$ and~$s$,
The exact sequences of~\ref{sbsc-fil0-str}
form a commutative diagram
\bel{diag-fil0-structure}
\begin{split}
\xymatrix{
0\ar[r] & H^q_{p^s}(\kk')\ar[r] & {\fil}_0 H^q_{p^s}(K')\ar[r]
  & H^{q-1}_{p^s}(\kk')\ar[r] & 0 \\
0\ar[r] & H^q_{p^s}(\kk)\ar[r]\ar[u] & {\fil}_0 H^q_{p^s}(K)\ar[r]\ar[u]
  & H^{q-1}_{p^s}(\kk)\ar[r]\ar[u]^{e} & 0,
}
\end{split}
\end{equation}
where the left two vertical maps are the canonical maps and the rightmost
vertical map is the canonical map multiplied by~$e$.
\end{prop}
\begin{proof}
Apply~\ref{sbsc-fil0-str}.
\end{proof}

\begin{prop}
\lb{prop-tame-action}
Suppose the extension~$A'/A$ is tame and generically
Galois with group~$G$.  Then its inertia group~$G_0$ acts trivially
on~$\fil_0 H^q_{p^s}(K')$.
\end{prop}
\begin{proof}
Assume, as we may, that~$A'/A$ is residually trivial, and
consider~(\ref{diag-fil0-structure}).
Since~$A'/A$ is tame and~$\kk'=\kk$, the outer
vertical morphisms are isomorphisms.  Therefore the inner vertical map
is, too.  Since~$G$ acts trivially on~$\fil_0 H^q_{p^s}(K)$, it acts
trivially on~$\fil_0 H^q_{p^s}(K')$.
\end{proof}

\section{The proof: equal characteristic}
\lb{sec-equal-kato}
\newcommand{\Hone}{{{H}^1(K,\qq_p/\zz_p)}}
\newcommand{\HH}{H}

In this section, we use Matsuda's refinement~\cite{Matsuda:Swan} of
Kato's refined Swan conductor to prove the the theorem in
equal characteristic.  Again, let us first recall the basics.

Assume~$A$ is of equal characteristic.

\subsection{}
For any positive integer~$s$,
let~$W_s(K)$ denote the group of Witt vectors of~$K$ of length~$s$.
The Verschiebung maps form an inductive system of abelian groups
\[
W_1(K)\longlabelmap{V} W_2(K)\longlabelmap{V} W_3(K)\longlabelmap{V}\cdots
\]
As in Fontaine \cite{Fontaine:Co-Witt}, let~$\CW(K)$ be its colimit.
For example, if~$\ff_p$ denotes the finite field of~$p$ elements,
then~$\CW(\ff_p)=\qq_p/\zz_p$.  The maps~$W_s(K)\to\prod_{-\nn}K$ defined by
\[
(a_{-s+1},\dots,a_0)\mapsto (\dots,0,a_{-s+1},\dots,a_0)
\]
induce a bijection between~$\CW(K)$ and the set of
elements~$(\dots, a_{-1}, a_0)$ such that~$a_{-i}=0$ for
sufficiently large~$i$.  I will typically use this identification without
comment.

\subsection{}
The Frobenius endomorphisms of the groups~$W_s(K)$ extend to an endomorphism
of~$\CW(K)$.  Call it~$F$.  For any separable closure~$K\sep$
of~$K$, we have an Artin-Schreier sequence of~$\Gal(K\sep/K)$-modules
\[
0\longmap\CW(\ff_p)\longmap\CW(K\sep)\longlabelmap{F-1}\CW(K\sep)\longmap 0.
\]
It is easy to show that this sequence induces a surjection
\[
\AS:\CW(K)\longmap\HH^1(K,\CW(\ff_p))=\HH^1(K,\qq_p/\zz_p)
\]
with kernel~$(F-1)\cdot\CW(K)$.

\subsection{}
Let~$\varphi:\CW(K)\rightarrow\Omega^1_K$ denote the homomorphism 
defined by
\[
(\dots,a_{-1},a_0) \longmapsto -\sum_{i\in\nn}a_{-i}^{p^i-1}da_{-i}.
\]
It appears that this map was first considered by Serre
\cite{Serre:Mexico}; it has a nice interpretation in terms of 
the de\! Rham-Witt complex \cite[I 3.12]{Illusie:dRW}.

\subsection{}
\lb{sbsc-Matsuda-filt}
Define the following filtrations indexed by non-negative
integers~$n$:
\[\begin{array}{lll}
\fillog_n \CW(K) & = &\{(\dots,a_{-1},a_0)\mid\forall i\in\nn\ p^{i}v_A(a_{-i})\geq-n \},\\
\fillog_n \Omega^1_K & = &\m_A^{-n}\cdot \dlog(K^*), \ \mbox{and} \\
\fillog_n \HH^1(K,\qq_p/\zz_p) & = & \AS(\fillog_n \CW(K)).
\end{array}\]
The filtration on fixed length Witt vectors was
apparently first considered by Schmid \cite{Schmid:Filt} in the
residually perfect case and (independently) Brylinski
\cite{Brylinski:Ram} in the residually imperfect case.  It is
immediate that~$\AS$ and~$\varphi$ preserve these filtrations.  It
is easy to check that~$\varphi$ does not factor
through~$\AS$ but that~$\grs\varphi$ does factor through~$\grs\AS$.
Matsuda
remarked~\cite[3.2.2]{Matsuda:Swan} that we have even more:

\begin{prop}
\label{matsuda-prop}
Let~$n$ be a non-negative integer.
Then~$\fillog_n\varphi/\fillog_{[n/p]}\varphi$
factors through~$\fillog_n\AS/\fillog_{[n/p]}\AS$.
\end{prop}

(Here,~$[x]$ is the greatest integer that is at most~$x$.)

\subsection{}
Denote by~$\phi_n$ the resulting homomorphism
\[
\fillog_n\Hone/\fillog_{[n/p]}\Hone\to\fillog_n\Omega^1_K/\fillog_{[n/p]}\Omega^1_K
\]
and by~$\grs\phi$ the induced map~$\grs\Hone\to\grs\Omega^1_K$ of
associated graded modules.  The following theorem of Kato's \cite[2.5, 3.2,
3.7]{Kato:Imp} shows how to use~$\grs\phi$ to compute the refined
Swan conductor of a Galois character given a representation of it as a Witt
vector. 

\begin{thm}
\lb{thm-eq-char-rsw}
For every positive integer~$n$, we have
\[
\fillog_n\HH^1(K,\qq/\zz)\cap\Hone = \fillog_n\Hone,
\]
where~$\Hone$ is viewed as the~$p^{\infty}$-torsion subgroup
of~$\HH^1(K,\qq/\zz)$.  Furthermore, under the natural
identification~$\gr_n\Omega^1_K=\m_A^{-n}\otimes_A\Omlog{1}{\kk}$, the
restriction of~$\kappa_{n}$
to~$\gr_n\Hone\subseteq\gr_n\HH^1(K,\qq/\zz)$ coincides
with~$\gr_n\phi$.
\end{thm}

\begin{prop}
\label{prop-compat}
For non-negative integers~$m\leq n$, the diagram
\[\xymatrix{
{\fillog}_n\Hone\ar[r]^{\phi_n} & {\fillog}_n\Omega^1_K/\fillog_{[n/p]}\Omega^1_K \\
{\fillog}_m\Hone\ar[u]\ar[r]^{\phi_m} & {\fillog}_m\Omega^1_K/\fillog_{[m/p]}\Omega^1_K\ar[u]
}\]
commutes.
\end{prop}
\begin{proof} Clear.
%The maps~$\fillog_m\varphi$ and~$\fillog_n\varphi$ agree
%on~$\fillog_m \CW(K)$ by definition.  Since the maps in question are the
%induced factors of these maps modulo~$\fillog_{[n/p]}\Omega^1_K$, the
%result follows.
\end{proof}

\begin{prop}
\label{prop-func}
Let~$A'$ be an extension of~$A$ of ramification index~$e$, and let~$K'$
denote its fraction field.  If~$\chi\in\HH^1(K,\qq_p/\zz_p)$, then
\[
\phi_{\swk(\chi)}(\chi)|_{A'} = \phi_{\swk(\chi|_{A'})}(\chi|_{A'})\Mod \fillog_{[e\,\swk(\chi)/p]}\Omega^1_{K'}.
\]
\end{prop}

\begin{proof}
%Since~$\swk(\chi|_{A'})\leq e\,\swk(\chi)$, this
%follows immediately from~\ref{prop-compat}.
By~\ref{prop-compat}.
\end{proof}

%(Compare with~\ref{prop-rsw-natural}.)

%This proposition is a further refinement of proposition
%\ref{lambda-natural} and is one reason why we have better results in
%equal characteristic.  Proposition \ref{lambda-natural} reduces to
%$0=0$ unless~$\swk(\chi|_{K'})=e\,\swk(\chi)$, when the ramification
%of $\chi$ is unchanged.  On the other hand, there is often enough
%room between~$\fillog_{en}$ and~$\fillog_{[en/p]}$ to see
%the refined Swan conductor of~$\chi|_{K'}$ in~$\phi_{\swk(\chi)}(\chi)|_{K'}$.

%This proposition is one reason why we have better results in equal
%characteristic.  In mixed characteristic, we have only proposition
%\ref{lambda-natural}, which reduces to $0=0$
%unless~$\swk(\chi|_{K'})=e\,\swk(\chi)$.  On the other hand, there is
%often enough room between~$\fillog_{en}\Omega^1_{K'}$
%and~$\fillog_{[en/p]}\Omega^1_{K'}$ to see the refined Swan conductor
%of~$\chi|_{K'}$ in~$\phi_{\swk(\chi)}(\chi)|_{K'}$.

\subsection{}
Finally, for any non-negative integer~$n$,
put~$\filreg{n}\Omega^1_K=\m_A^{-n}\cdot \Omega^1_A\subset\Omega^1_K$.
This filtration measures the order of the pole in the usual sense,
whereas~$\fil_{\scriptdot}\Omega^1_K$ measures it in the
logarithmic sense.  The two filtrations are intertwined:
\[
\cdots\subseteq\fillog_{n}\Omega^1_K\subseteq\filreg{n+1}\Omega^1_K\subseteq\fillog_{n+1}\Omega^1_K\subseteq\cdots.
\]

Matsuda~\cite[3.1]{Matsuda:Swan} has given non-logarithmic variants of
the other filtrations in~\ref{sbsc-Matsuda-filt},
but we will not need them here.  (Note,
however, that our indexing of~$\filreg{\scriptdot}\Omega^1_K$ differs from
Matsuda's by one.)

\begin{prop}
\lb{prop-eq-char-ark}
Let~$\chi$ be a class in~$\Hone$.  Then~$\ark(\chi)$ is
the smallest integer~$n$ satisfying~$\chi\in\filreg{n}\Hone$.
\end{prop}
\begin{proof}
By~\ref{thm-eq-char-rsw}.
\end{proof}

\begin{lemma}
\lb{lemma-delta-inj2}
For~$n\geq 1$, the natural map~$\gr'_n\Omega^1_K
\longmap\gr'_n\Omega^1_{K\gc}$ is injective.
\end{lemma}
\begin{proof}
Since we have~$\gr'_n\Omega^1_K = \Omega_A\otimes_A \p_A^n/\p_A^{n-1}$,
it is enough to show the map
\[
\kk\otimes_A \Omega^1_A \longmap\kk\otimes_A\Omega^1_{A\gc}
\]
is injective.

Let~$T$ be a lift to~$A$ of a $p$-basis for~$\kk$, and let~$\pi$
be a uniformizer of~$A$.  Then the set~$dT\cup \{d\pi\}$ is a basis for
the~$\kk$-module~$\kk\otimes_A\Omega^1_A$.  To show
injectivity, it is enough to check that the image of~$dT\cup \{d\pi\}$
is~$\kk$-linearly independent.  But~$\Omega^1_{A\gc}=A\gc d\pi$,
and so it is enough to show, in the notation of~\ref{thm-expl},
that~$\{1\}\cup\{u_{t,1}\mid t\in T\}$ is~$\kk$-linearly
independent in~$k\gc$.  This follows from~\ref{thm-expl}.
\end{proof}

We can now prove the the theorem in equal characteristic.

\begin{proof}
First suppose~$\chi$ is in~$H^1(K,\qq_p/\zz_p)$. Put~$n=\swk(\chi)$
and~$m=\swk(\chi|_{K\gc})$. 
Then~$\ark(\chi)$ is either~$n$ or~$n+1$.  We will treat
these two subcases separately.  If~$\ark(\chi)=n+1$,
then~$\kappa_n(\chi)\notin\m_A^{-n}\otimes\Omega^1_{\kk}$.  The
naturality~(\ref{prop-rsw-natural}) of the maps~$\kappa_{\scriptdot}$
therefore implies~$\kappa_n(\chi|_{A\gc})\notin\m_{A\gc}^{-n}\otimes
\Omega^1_{k\gc}$ and, so,~$\ark(\chi|_{A\gc})=n+1$.

Now consider the second subcase, when~$\ark(\chi)$ is~$n$. 
By~\ref{prop-p-div-sw}, we
have~$n\geq 2$; so, for~$n=2$, it is enough to show~$\chi|_{A\gc}$
is not tame.  This follows from~\ref{thm-Galois-surj}.
If~$n \geq 3$, we have~$[n/p]\leq n-2$ and, hence,
\[
\fil_{[n/p]}\Omega^1_{K\gc} \subseteq \fil_{n-2}\Omega^1_{K\gc}
  \subseteq \filreg{n-1}\Omega^1_{K\gc}.
\]
Then, by~\ref{prop-func} and~\ref{lemma-delta-inj2}, we have
\[
\phi_m(\chi|_{A\gc}) \equiv \phi_n(\chi)|_{A\gc} \not\equiv 0
\mod \filreg{n-1}\Omega^1_{K\gc}.
\]
Therefore,~$\chi\not\in\filreg{n-1}\Omega^1_{K\gc}$,
and so~\ref{prop-eq-char-ark} implies~$\ark(\chi)=n$.

Now let~$\chi$ be an arbitrary class in~$H^1(K)$.  If~$\chi$ is tame,
the result follows from~\ref{thm-Galois-surj}.  If~$\chi$ is wild,
then write~$\chi=\chi'+\chi''$, where~$\chi'$ is
in~$H^1(K,\qq_p/\zz_p)$ and~$\chi''$ is tame.  Since~$\chi''$ is
tame,~$\chi$ and~$\chi'$ have the the same refined Swan conductor and,
hence, the same Kato-Artin conductor.  Similarly,~$\chi'|_{K\gc}$ is
wild (again by~\ref{thm-Galois-surj}), and so~$\chi|_{A\gc}$
and~$\chi'|_{A\gc}$ have the same Kato-Artin conductor.  The work
above then implies~$\ark(\chi')=\ark(\chi'|_{A\gc})$, and this completes
the proof.
\end{proof}

\section{Some lemmas}

The purpose of this section is to prove some lemmas needed in the proof
of the the theorem in mixed characteristic.  All the results in this section
are, however, still valid in equal characteristic.
Let~$\pi$ be a uniformizer of~$A$.

\subsection{}
\lb{sbsc-K-thy}
Let~$U^{i}\KM{2}(K)$ (for~$i\geq 1$) denote the
subgroup of~$\KM{2}(K)$ (see~\ref{sbsc-Kato-notation}) generated by the set~$\{U^{i}_A, K^*\}$.  This
filtration satisfies~\cite[4.1]{Bloch-Kato:p-adic}
\bel{equ-K-fil}
\{U^i_A,U^j_A\} \subseteq U^{i+j}\KM{2}(K).
\end{equation}

\begin{lemma}
\lb{lemma-symbol-identity}
Let~$x$ and~$y$ be non-zero elements of~$\m_A$.  Then
\[
\{1+ x, 1+y \} \equiv  \{1+xy,-y\} \mod U^{v_A(xy)+1}\KM{2}(K).
\]
\end{lemma}
\begin{proof}
We have
\begin{align*}
\{1+x,1+y\} &= \{-y(1+x),1+y\} \\
            &= -\{-y(1+x),1+xy(1+y)^{-1}\} \\
            &\equiv -\{-y(1+x),1+xy\} \mod U^{v_A(xy)+1}\KM{2}(K) \\
            &\equiv \{1+xy,-y\} \mod U^{v_A(xy)+1}\KM{2}(K).
\end{align*}
\end{proof}

\begin{lemma}
\lb{lemma-submain}
Suppose~$A$ is residually perfect.
Let~$x\neq 0$ be an element of~$\m_A$, and
let~$z\in A^*$ be such that the
element~$z'=z-s_A(\bar{z})$ is non-zero. Then
\[
\{1+ x, z \} \equiv  v_A(z')\{1+xz's_A(\bar{z}^{-1}),\pi\}
   \mod U^{v_A(xz')+1}\KM{2}(K) + D\KM{2}(K),
\]
where~$D\KM{2}(K)$ is the infinitely $p$-divisible subgroup of~$\KM{2}(K)$.
\end{lemma}

\begin{proof}
The defining property of the lift~$s_A(\bar{z})$ of~$\bar{z}$ is that it is
infinitely $p$-divisible in~$K^*$.  Therefore, it suffices to
assume~$\bar{z}=1$.  By~\ref{lemma-symbol-identity}, we have
\[
\{1+ x, 1+z' \} \equiv  v_A(z')\{1+xz',\pi\} + \{1+xz',-z'/\pi^{v_A(z')}\}
   \mod U^{v_A(xz')+1}\KM{2}(K).
\]
Since~$-z'/\pi^m$ is in~$A^* = s_A(\kk^*) U^1_A$, we also have
\[
\{1+xz',-z'/\pi^m\}\in U^{v_A(xz')+1}\KM{2}(K) + D\KM{2}(K).
\]
\end{proof}

%\begin{lemma}
%\lb{lemma-fil-shift}
%Let~$m,n$ be non-negative integers with~$n\geq m$,
%and let~$\chi$ be an element of~$\fil_nH^q(K)$.
%If~$x\in U^m_A$, then~$\{\chi,x\}\in\fil_{n-m}H^{q+1}(K)$.
%\end{lemma}
%\begin{proof} 
%Write~$x=1+y$.  Then for any henselian extension~$A'/A$ and any
%element~$z$ in~$\m_{A}^{n-m}\m_{A'}$, we have
%by~\ref{lemma-symbol-identity}
%\begin{eqnarray*}
%\{\chi,x,1+z\} & = & \{\chi,1+y,1+z\} \\
%               & = & \{\chi,1+yz,-z\} \\
%               & = & 0, 
%\end{eqnarray*}
%and therefore~$\{\chi,x\}$ is in~$\fil_{n-m}H^{q+1}(K)$.
%\end{proof}

\begin{prop}
\lb{prop-coh-p-div}
If~$F$ is a field of characteristic~$p$, the $p^\infty$-torsion
subgroup of~$H^q(F)$ is $p$-divisible.
\end{prop}
\begin{proof}
%This is an immediate consequence of 
By the surjectivity of the higher
Artin-Schreier maps~$\AS_{q}$ and the natural
maps~$\Omega^{\scriptdot}_{W_s(F)}\to W_s\Omega^{\scriptdot}_F$
and~$W_{s+1}(F)\to W_s(F)$.
\end{proof}

\begin{lemma}
\lb{lemma-not-in-fil0}
Let~$\chi\in H^1(K)$ be a class such that~$\ark(\chi)=\swk(\chi)\neq 0$.
Then, in the notation of~\ref{sbsc-tilde}, the element
%\[
$\{\chi|_{\tilde{A}},1+\pi^{\swk(\chi)-1}T\}$
%\]
is in~$\fil_1H^2_p(\tilde{K}) + H^2_{p^2}(\kkt)$ but is
not in~$\fil_0 H^2(\tilde{K})$.
\end{lemma}
\begin{proof}
Write~$n=\swk(\chi)$.  Then by~\ref{prop-p-div-sw}, we have~$n>1$.
Let~$\psi = \{\chi|_{\tilde{A}},1+\pi^{n-1}T\}$.
Since~$n>1$, we have
\[
(1+\pi^{n-1}T)^{p^2}\in U^{n+1}_{\tilde{A}}\ \ \text{and}\ \ (1+\pi^{n-1}T)^p\in U^{n}_{\tilde{A}}.
\]
Applying~(\ref{equ-K-fil}),
we get~$\psi\in\fil_1 H^2_{p^2}(\tilde{K})$ and~$p\psi\in\fil_0
H^2_p(\tilde{K})$.  
%It follows from~(\ref{equ-K-fil})
%that~$\psi$ is contained in~$\fil_1 H^2(\tilde{K})$.

Let~$A'$ be a henselian extension of~$A$ whose
residue field~$k'$ is perfect and has the property that~$k$ is separably
closed in it.  (Take~$A'=A\gc$, for example.)
Since~$\ark(\chi)=\swk(\chi)$ and since~$k'$ is
perfect,~\ref{prop-rsw-natural} implies~$\swk(\chi|_{A'})\leq n-1$.
Therefore,~$p\psi|_{\widetilde{A'}}$ is zero.  Since~$k$ is
separably closed in~$k'$, the field~$\tilde{k}=k(T)$ is separably
%\marginpar{why?}
closed in~$\tilde{k'}=k'(T)$.
The natural map~$H^1_p(\tilde{k}) \to H^1_p(\tilde{k'})$ is therefore an injection.
So, by chasing diagram~(\ref{diag-fil0-structure}) applied to the
extension~$\tilde{A'}/\tilde{A}$, we conclude~$p\psi \in H^2_p(\tilde{k})$.
%Therefore, the rightmost vertical map
%of~(\ref{diag-fil0-structure}), applied to the extension~$A'/A$, is an
%injection.  So, by chasing the diagram, we
%see~$p\psi \in H^2_p(\tilde{k})$.
By~\ref{prop-coh-p-div}, there is
a class~$\psi''\in H^2_{p^2}(\kkt)$ such that~$p\psi''=p\psi$.
Then,~$\psi-\psi''$ is in~$\fil_1 H^2_p(\tilde{K})$.

On the other hand, putting~$\eta=\kappa_n(\chi)$, we have
\begin{equation*}
\begin{aligned}
\{\psi,1-\pi\} &= \{\chi|_{\tilde{A}},1+\pi^{n-1}T,1-\pi\} \\
   &= \{\chi|_{\tilde{A}},1-\pi^{n}T,\pi\} \\
   &= \{\lambda_{\tilde{A}}(-\eta \pi^n T), \pi\}.
\end{aligned}
\quad
\begin{gathered}
\\
\text{by~\ref{lemma-symbol-identity}} \\
\text{by~\ref{thm-kato-thy}.}
\end{gathered}
\end{equation*}
But, by assumption,~$\eta$ is in~$\m_A^{-n}\otimes{\Omega}^1_{\kk}$,
and so~$\lambda_{\tilde{A}}(-\eta\pi^n T)$ equals~$\AS_1(-\eta\pi^n
T)$, which is non-zero~\cite[3.8]{Kato:Imp}.
%Thus,~$\{\lambda_{\tilde{A}}(-\eta \pi^n T), \pi\}$ is non-zero
Because of this,~$\{\psi,1-\pi\}$ is non-zero
(\ref{sbsc-fil0-str}) and, hence,~$\psi$ is not in~$\fil_0
H^2(\tilde{K})$.
\end{proof}

\section{A diagram}
\lb{sec-lemmas}
\newcommand{\abar}{{\kk}}

Assume~$A$ is of mixed characteristic.
Let~$\KMb{2}(K)$ denote~$\KM{2}(K)/p\KM{2}(K)$, and
let~$U^{\scriptdot}\KMb{2}(K)$ denote the image of the
filtration~$U^{\scriptdot}\KM{2}(K)$. 

The primary purpose of this section is to prove a certain diagram
(\ref{prop-main}) commutes.  This will allow us to understand in terms
of symbols how some classes in~$H^2_p(K)$ change when pulled back
to~$A\gc$.  The relation between symbols and cohomology classes is
provided by a theorem from
Bloch-Kato~\cite[5.12]{Bloch-Kato:p-adic} and
Kato~\cite[4.1(6)]{Kato:Imp}:

\begin{thm}
If~$A$ contains the group~$\mu_p$ of all $p$-th roots of unity,
then~$h_K$ (of~\ref{sbsc-Kato-notation}) induces isomorphisms
\[
\KMb{2}(K) \longisomap H^2_p(K)\otimes\cech{\mu}_p
\quad\text{and}\quad
U^{\eh-n}\KMb{2}(K)\otimes\cech{\mu}_p \longisomap {\fil}_n H^2_p(K),
\]
where~$\eh=e_{A/\zz_p}p(p-1)^{-1}$ and~$0\leq n < \eh$.
\end{thm}

\subsection{}
Fix an extension~$A_\mu$ of~$A$ that is generically generated by a
primitive~$p$-th root of unity.  We will assume thoughout this section
that the extension~$A_\mu/A$ is residually trivial.
%For the rest of this section, we assume only that~$A$ admits a
%residually trivial extension that contains a primitive~$p$-th root of unity.
Let~$A'$ be a henselian residually perfect extension of~$A$ of
ramification index one.  Because we will use this only when~$A'=A\gc$,
the reader is free to assume it (even though it does not
simplify anything).

Aso fix the following notation:
% $A_\mu/A$ is an extension generically
%generated by a primitive $p$-th root of unity,
$A'_\mu$ is $A'\otimes_A A_\mu$ (a
residually perfect henselian discrete valuation
ring);~$K',K_\mu,K'_\mu$ are the fraction fields of~$A',A_\mu,A'_\mu$,
and~$k'$ is the residue field of~$A'$; $\Gamma$
is~$\Gal(K_\mu/K)$, and~$m$ is the degree of~$K_\mu/K$; $\mu_p$ is
the group of all~$p$-th roots of unity in~$K_\mu$, and~$\cech{\mu}_p$
is~$\Hom(\mu_p,\zz/p)$; $\zeta$ is some
non-trivial element of~$\mu_p$, and~$\cech{\zeta}$ is the element
in~$\cech{\mu}_p$ such that~$\cech{\zeta}(\zeta)=1$.  (None of the
constructions below will depend on the choice of~$\zeta$.)

\begin{prop}
\lb{prop-gamma}
There is a map~$\gamma_{A'/A}$ that makes the following diagram commute:
\[\xymatrix{
{\fil}_mH^2_p(K_\mu) \ar[d]\ar@{>>}[r]
  &  {\gr}_mH^2_p(K_\mu) \ar@{-->}[ddd]^-{\gamma_{A'/A}}
  \\
H^2_p(K'_\mu) \\
{\fil}_0 H^2_p(K'_\mu) \ar@2{-}[u]\ar[d]%^-{\simeq}
  \\
H^1_p(k') \ar@{>>}[r]
  &  \coker{H^1_p(k)\to H^1_p(k')}.
}\]
\end{prop}

(Note that~$\fil_0 H^2_p(K'_\mu)=H^2_p(K'_\mu)$
since~$A'_\mu$ is residually perfect~\cite[6.1]{Kato:Imp}.)

\begin{proof}
Let~$\chi$ be a class in~$\fil_{m-1}H^2_p(K_\mu)$, and
put~$\chi'=m^{-1}\sum_{\sigma\in\Gamma}\sigma(\chi)$.
%\[ 
%m^{-1}\sum_{\sigma\in\Gamma}\sigma(\chi).
%\]
It is clear that~$\chi'$ lies in~$\fil_{m-1}H^2_p(K_\mu)^\Gamma$, which
by~\ref{prop-G-inv}, agrees with~$\fil_0 H^2_p(K)$.
On the other hand,~\ref{prop-tame-action} implies~$\Gamma$ acts trivially
on~$\fil_0 H^2_p(K'_\mu)=H^2_p(K'_\mu)$, and so~$\chi$ and~$\chi'$ have
the same image in~$H^2_p(K'_\mu)$.  The image of~$\chi$ in~$H^1_p(k')$
is therefore contained in the image of~$H^1(k)$.
\end{proof}

%\pagebreak
\begin{prop}
\lb{prop-der}
The composite map
\[
A \to \m_{A'}/(\m_{A'}^2+\m_A) = (\kk'/\kk) \otimes_A\m_A
\longmap (\kk'/\proot{\kk}) \otimes_A\m_A,
\]
where the leftmost map sends~$x$ to the class
of~$x - s_{A'}(\bar{x})$, is a derivation that vanishes on~$\m_A$.
\end{prop}
\begin{proof}
It is immediate that it vanishes on~$\m_A$ and a short
computation shows it satisfies the Leibniz rule.
To see it is additive, it suffices to show that for
all~$x,y\in\kk$, we have
\begin{equation*}
s_{A'}(x+y) \equiv s_{A'}(x) + s_{A'}(y)\mod\m_{A'}^2 + s_{A'}(\proot{k})\m_A.
\end{equation*}
Let~$W$ be the ring of Witt
vectors~\cite[0.1]{Illusie:dRW} with entries in~$k^{p^{-\infty}}$.  By the
definition of addition in~$W$, there is an element~$z\in\kk$ such
that
\[
(x,0,\dots) + (y,0,\dots) = (x+y,z,\dots),\text{\ \ i.e.,}
\]
\[
s_{W}(x) + s_{W}(y) \equiv s_{W}(x+y) + ps_{W}(\proot{z}) \mod p^2.
\]
Since there is a map~$W\to A'$ that is compatible with multiplicative
sections, the congruence above holds.
\end{proof}

If~$p$ does not generate~$\m_A$, even the
map~$A\to(\kk'/\kk) \otimes_A\m_A$ is a derivation (in both
equal and mixed characteristic).

\begin{const}
{\em Morphisms~$\del_{A'/A},\nu,$ and~$\theta$}
%{\em Morphisms
%$\del_{A'/A}:\m_A^{-1}\otimes_A\Omega^1_\kk \longmap \kk'/\proot{\kk},
%\nu:{\gr}^{\eh-m}\KMb{2}(K_\mu)\otimes\cech{\mu}_p \longmap
%{\m}_{A}^{-1}\otimes\Omega^1_{\kk},$ and~$\theta: k' \longmap U^{\eh}\KMb{2}(K'_\mu)\otimes\cech{\mu}_p$}

Let~$\del_{A'/A}:\m_A^{-1}\otimes_A\Omega^1_\kk \longmap \kk'/\proot{\kk}$
denote the $\kk$-linear homomorphism induced by the derivation
in~\ref{prop-der}.  Put
\[
\eh = e_{A_\mu/\zz_p}p(p-1)^{-1} = v_{A_\mu}\bigl(p(\zeta-1)\bigr),
\]
and let (in the notation of~\ref{sbsc-K-thy})
\[
\nu:{\gr}^{\eh-m}\KMb{2}(K_\mu)\otimes\cech{\mu}_p \longmap
{\m}_{A}^{-1}\otimes\Omega^1_{\kk} 
\]
denote the map determined by
\[
\{1-p(\zeta-1)x,y\}\otimes\cech{\zeta}  \mapsto
mx\otimes\dlog(\bar{y}),
\]
where~$x\in\m_A^{-1}A_\mu$,~$y\in A_\mu$.  Bloch and Kato show
%\marginpar{Check.}
it is (well-defined and) an isomorphism~\cite[4.3, 5.2]{Bloch-Kato:p-adic}.
Because~$U^{\eh+1}_{A_\mu}\subseteq(K_\mu^*)^p$,
we have~$U^{\eh}\KMb{2}(K'_\mu)=\gr^{\eh}\KMb{2}(K'_\mu)$.
We can therefore define a map
\[
\theta: k' \longmap U^{\eh}\KMb{2}(K'_\mu)\otimes\cech{\mu}_p
\]
by~$x\mapsto \{1-p(\zeta-1)x,\pi_\mu\}$, where~$\pi_\mu$ is
a uniformizer of~$A_\mu$.  (The map is independent of the choice.)
\end{const}

\begin{prop}
\lb{prop-main}
The following diagram commutes:
\[
\begin{xy} <0pt,0pt>;<17.5mm,0mm>:<0mm,-7mm>::
  (2,0) *+{k'/k^{p^{-1}}}="a",
  (5,0) *+{\coker{H^1_p(k^{p^{-1}})\to H^1_p(k')}}="b",
  (5,3) *+{\coker{H^1_p(k)\to H^1_p(k')}}="d",
  (0,4) *+{{\m}_A^{-1}\otimes_A\Omega^1_{k}}="e",
  (2,6) *+{k'}="f",
  (5,6) *+{H^1_p(k')}="g",
  (0,7) *+{{\gr}^{\eh-m}\KMb{2}(K_\mu)\otimes\cech{\mu}_p}="h",
  (3,7) *+{{\gr}_m H^2_p(K_\mu)}="i",
  (1,8) *+{U^{\eh}\KMb{2}(K'_\mu)\otimes\cech{\mu}_p}="j",
  (4,8) *+{{\fil}_0 H^2_p(K'_\mu)}="k",
  (0,10) *+{U^{\eh-m}\KMb{2}(K_\mu)\otimes\cech{\mu}_p}="l",
  (3,10) *+{{\fil}_m H^2_p(K_\mu).}="m",
  \ar @{->>}"a";"b" ^-{\AS_1}
  \ar "d";"b" ^-{\simeq}
  \ar "e";"a" ^-{{\del}_{A'/A}}
  \ar @{->>}"f";"a"
  \ar @{->>}"f";"g" ^(.37){{\AS}_1} |!{"i";"d"}\hole
  \ar @{->>}"f";"j" _(.3){\theta} |!{"h";"i"}\hole
  \ar @{->>}"g";"d"
  \ar "h";"e" ^-{\nu}_-{\simeq}
  \ar "h";"i" ^(.65){\simeq}_(.65){h_{K_\mu}}

  \ar "i";"d" _-{\gamma_{A'/A}}
  \ar "j";"k" ^-{\simeq}_-{h_{K'_\mu}} |!{"m";"i"}\hole
  \ar "k";"g" ^-{\simeq}
  \ar "l";"j"
  \ar @{->>}"l";"h"
  \ar "l";"m" ^-{\simeq}_-{h_{K_\mu}}
  \ar @{->>}"m";"i"
  \ar "m";"k"
\end{xy}
\]
\end{prop}
\begin{proof}
The commutativity of the rear lower face follows from the splittings
of~\ref{sbsc-fil0-str} and the usual compatibility between
Artin-Schreier theory and Kummer theory. It is clear the other three
faces for which it makes sense to ask the question commute.
Therefore, it only remains to check that the perimeter commutes.

Because~$\nu$ is an isomorphism, it is enough to consider elements
\[
x=\nu^{-1}(\pi^{-1}\otimes dz) = \{1-p(\zeta-1)\pi^{-1}z,z\}\otimes\cech{\zeta}
  \in {U^{\eh-m}\KMb{2}(K_\mu)\otimes\cech{\mu}_p}.
\]
where~$\pi$ is a uniformizer of~$A$
and~$z\in A^*$.  Write~$z|_{A'} = s_{A'}(\bar{z}) + \pi y$,
where~$y\in A'$.  Then, letting~$[x]$ denote the graded class of~$x$, we have
\[
\del_{A'/A}\circ\nu([x]) = \bar{y} \bmod \proot{\kk}.
\]
On the other hand, by~\ref{lemma-submain}, we have
\[
x|_{K'_\mu} = \{1-p(\zeta-1)yzs(\bar{z})^{-1},\pi_\mu\} \otimes \cech{\zeta}
               = \theta(\bar{y}).
\]
Since the front lower and rear faces commute, the proof is complete.
\end{proof}

\section{The proof: mixed characteristic}
\lb{sec-mixed-kato}

Assume in this section that~$A$ is of mixed characteristic, and
let~$A_\mu$ be an extension of~$A$ that is generically generated by a
primitive~$p$-th root of unity.  (The extension~$A_\mu/A$ is no longer
assumed to be residually trivial.)

\begin{prop}
\lb{prop-del-inj}
The composite map
\[\xymatrix{
{\m}_A^{-1}\otimes_A\Omega^1_{k} \ar[r]^-{\del_{A\gc/A}}
  & k\gc/\proot{k} \ar[r]^-{\AS_1}
  & \coker{H^1_p(\proot{k}) \to H^1_p(k\gc)}
}\] 
is injective.
\end{prop}
(Compare with~\ref{lemma-delta-inj2}.)
\begin{proof}
Let~$T$ be a $p$-basis for~$\abar$.  Then~$dT$ is a basis for~$\Omega^1_\abar$.
%Let~$\pi$ be a uniformizer of~$A$.
Let~$\eta$ be an element of the kernel of~$\AS_1\circ \del_{A\gc/A}$ and write
\[
\eta =  \pi^{-1}\otimes \sum_{t\in T}a_t dt,\ a_t\in k
\]
where~$\pi$ is a uniformizer of~$A$ and~$a_t$ is zero
for all but finitely many~$t\in T$.
By~\ref{thm-expl}, we have
\[
\del_{A\gc/A}(\eta) = \sum_{t\in T} a_t u_{t,1}\bmod\proot{\abar}.
\]
%the image of~$\eta$ in~$\abar\gc/\proot{\abar}$
%is~$\sum_t a_t u_{t,1}\bmod\proot{\abar}$.
Since the image of this in
\[
\coker{H^1_p(\proot{\abar}) \to H^1_p(\abar\gc)}
\]
is assumed to be
zero, there are elements~$x\in \abar\gc$ and~$y\in \abar$ such that
\[
x^p-x = y + \sum_{t\in T} a_t u_{t,1}.
\]

Suppose for a contradiction that there is an element~$s\in T$ such
that~$a_s\neq 0$. Put
\[
F=\abar(u_{t,1}\mid t\neq s)^{p^{-\infty}}.
\]
Then~\ref{thm-expl} implies~$F(u_{s,1})$ is separably closed in~$\abar\gc$.
Therefore, we have~$x\in F(u_{s,1})$ and~$x^p - x = y' + a_s u_{s,1}$,
where~$y'=y+\sum_{t\neq s}a_tu_{t,1}\in F$.
But valuation considerations in the completion~$F((u_{s,1}^{-1}))$ 
show this is impossible.
\end{proof}

\begin{cor}
\lb{cor-gamma-inj}
If~$A_\mu/A$ is residually trivial, the map~$\gamma_{A\gc/A}$ is injective.
\end{cor}
\begin{proof}
By~\ref{prop-main} and~\ref{prop-del-inj}.
\end{proof}

\begin{prop}
\lb{prop-gc-H-inj}
Let~$s\geq 1$ be an integer.  Then the following sequences are exact:
\begin{equation}
\label{eq-gcHinj-1}
0  \longmap  H^2_{p^s}(k)  \longmap  \fil_0 H^2_{p^s}(K)
   \longmap  H^2_{p^s}(K\gc), 
\end{equation}
\begin{equation}
\label{eq-gcHinj-2}
0  \longmap  \hspace{1pt}H^2_{p}(k)\hspace{1pt}  \longmap  \hspace{2pt}\fil_1 H^2_{p}(K)\hspace{1pt}
   \longmap  H^2_{p}(K\gc).
\end{equation}
\end{prop}
\begin{proof}
Because~$k\gc$ is perfect, we have~$H^2_{p^s}(k\gc)=0$ (by, say, higher
Artin-Schreier theory~\cite[1.3]{Kato:Imp}).
Now,~\ref{thm-Galois-surj} implies the
map~$H^1_{p^s}(k)\to H^1_{p^s}(k\gc)$ is an injection.
%By~\ref{prop-fil0-structure}, we then have a map of exact sequences:
%\[\xymatrix{
%0\ar[r] & 0\ar[r]
%  & \fil_0H^2_{p^s}(K\gc)\ar[r]
%  & H^1_{p^s}(k\gc)\ar[r] & 0 \\
%0\ar[r] & H^2_{p^s}(k)\ar[r] \ar[r]\ar[u]
%  & \fil_0H^2_{p^s}(K)\ar[r]\ar[u]
%  & H^1_{p^s}(k)\ar[r]\ar[u] & 0.
%}\]
The exactness of~(\ref{eq-gcHinj-1})
then follows from~\ref{prop-fil0-structure}.

Let~$A_0$ be the maximal unramified subextension of~$A_\mu/A$.
Applying~\ref{cor-gamma-inj} to the extension~$A_\mu/A_0$,
we conclude that the map~$\gamma_{A\gc_0/A_0}$ is
injective, and therefore the kernel of the map
\[
\fil_m H^2_p(K_\mu) \longmap H^2_p(K\gc_0\otimes_{K_0}K_\mu)
\]
is contained in~$\fil_{m-1}H^2_p(K_\mu)$.
By~\ref{cor-tame-bc}, we have
\[
\fil_{m-1}H^2_p(K_\mu) \cap H^2_p(K) = \fil_0H^2_p(K).
\]
Because there is~\cite[2.1]{Borger:Cond-moduli} a map~$K\gc\to K\gc_0$,
we then have
\begin{align*}
\kernel{\fil_1H^2_p(K)\to H^2_p(K\gc)} & \subseteq
   \kernel{\fil_m H^2_p(K_\mu) \to H^2_p(K\gc_0\otimes_{K_0}K_\mu)}\cap H^2_p(K)\\
   & \subseteq \kernel{\fil_0H^2_p(K)\to H^2_p(K\gc)}, \\
   & = H^2_p(k),
\end{align*}
which proves the exactness of~(\ref{eq-gcHinj-2}).
\end{proof}

We can finally prove the theorem from the introduction
in mixed characteristic. 

\begin{proof} Let~$\pi\in A$ be a uniformizer.  Put~$n=\swk(\chi)$.
If~$n=0$, then~$\chi$ is tame, and the result follows immediately
from~\ref{thm-Galois-surj}.  Now assume~$n>0$.

By definition,~$\ark(\chi)$ is either~$n$ or~$n+1$.  If it
is~$n+1$, then~$\kappa_n(\chi)\notin\m_A^{-n}\otimes_A\Omega^1_{k}$.
Therefore,~\ref{prop-rsw-natural}
implies~$\kappa_n(\chi|_{A\gc})\notin\m_{A\gc}^{-n}\otimes_A
\Omlog{1}{k\gc}$, and so~$\ark(\chi|_{A\gc})=n+1$, as desired.

Now consider the case~$\ark(\chi)=n$.  Since~$A\gc$ is residually
perfect,~$\ark(\chi|_{A\gc})=n$ if and only
if~$\swk(\chi|_{A\gc})=n-1$.  By~\ref{prop-rsw-natural},
we know~$\swk(\chi|_{A\gc})$ is at most~$n-1$.  Put~$\psi =
\{\chi|_{\tilde{A}},1+\pi^{n-1}T\}\in\HH^2(\tilde{K})$.
Then~$\swk(\chi|_{A\gc})\geq n-1$ if and only
if~$\psi|_{\widetilde{A\gc}}$ is not zero.

Now we will construct a map~$\widetilde{A\gc}\longmap\tilde{A}\gc$.
(See figure~\ref{fig-mix}.)
Since~$\tilde{A}/A$ is residually separable,
%\marginpar{ref-prop-gc-func} %% Unstable reference
there is~\cite[2.1]{Borger:Cond-moduli} an
injection~$A\gc\to\tilde{A}\gc$.
Sending~$T\mapsto T$ yields another injection~$A\gc[T] \longmap \tilde{A}\gc$.
By the universal properties of localization and
henselization~\cite[VIII]{Raynaud:Hens}, this naturally induces a 
map~$\widetilde{A\gc}\longmap\tilde{A}\gc$.
It is therefore enough to show~$\psi|_{\tilde{A}\gc}\neq 0$.

\begin{figure}[htb] % htb
\[\xymatrix@=12pt{
& {\tilde{A}}\gc \\
&& {\widetilde{A\gc}}\ar@{.>}[ul] \\
{\tilde{A}}\ar[uur]\ar[urr] \\
&& A\gc\ar[uu] \ar[uuul]|!{[uu];[ull]}\hole\\
A \ar[uu]\ar[urr]
}\]
\caption{}
\label{fig-mix}
\end{figure}

%By~\ref{prop-p-div-sw}, we have~$n>1$.
By~\ref{lemma-not-in-fil0}, we can write~$\psi = \psi' + \psi''$,
where~$\psi'$ is in~$\fil_1H^2_p(\tilde{K})-\fil_0 H^2_p(\tilde{K})$
and~$\psi''$ is in~$H^2_{p^2}(\tilde{k})$.
We then have~$\psi|_{\tilde{A}\gc} = \psi'|_{\tilde{A}\gc} \neq 0$
by~\ref{prop-gc-H-inj}.
\end{proof}

\section{Comparison with Kato's conductor}
\lb{sec-cond}

Let~$\Lambda$ be a field whose characteristic is not~$p$, and fix an
injection from the torsion subgroup of~$\Lambda^*$ to~$\qq/\zz$.  In
each of the two results below,~$B$ is a finite extension of~$A$
that is generically Galois with
group~$G$, $\chi$ is a class in~$H^1(K,\qq/\zz)$, and
$\rho:G\to\Lambda^*$ is the corresponding homomorphism.
Let~$G_{\scriptdot}$ denote the lower ramification
filtration~\cite[IV Prop.\ 1]{Serre:CL} of~$G$, and
let~$\art(\rho)$ denote the the
Artin conductor~\cite[3.2]{Borger:Cond-moduli} of~$\rho$.
%% Unstable reference

\begin{prop}
\lb{prop-ark=arn}
If~$B/A$ is residually separable, we have
\[
\ark(\chi)=e_{B/A}^{-1}\sum_{i\geq 0}\card{G_i}\codim \Lambda^{G_i}
 = e_{B/A}^{-1}\sum_{\rho(G_i)\neq 1}\card{G_i}
\]
\end{prop}
\begin{proof}
If~$\rho$ is tame, the result is clear.  Now assume~$\rho$ is wild and
let~$n=\swk(\chi)$.  Then Kato
\hbox{\cite[6.8]{Kato:Imp}\cite[3.6(1),3.16]{Kato:Diff}}
shows~$\kappa_n(\chi)$ is not in~$\m_A^{-n}\otimes{\Omega}^1_{\kk}$
and that~$n$ agrees with the naive
Swan conductor~\cite[6.7.1]{Kato:Imp}.  And this, in turn, agrees
\cite[2.1]{Serre:Fact-locaux} with
\begin{equation} \label{naive-swan}
{e^{-1}_{B/A}}\sum_{i\geq 1}\card{G_i}\codim \Lambda^{G_i}.
\end{equation}
Since~$\ark(\chi)$ is one more than~$\swk(\chi)$ and since the
number~$\arnv{}(\rho)=\arnv{B}(\rho)$ is one more than
(\ref{naive-swan}), the equality of conductors follows.
\end{proof}

\begin{cor}
\lb{cor-ark=art}
$\art(\rho)=\ark(\chi)$.
\end{cor}

\begin{proof}
Since~$A\gc$ is residually perfect, %\marginpar{ref} %% Unstable reference
$\art(\rho)$ agrees~\cite[3.3]{Borger:Cond-moduli} 
with the sum in~\ref{prop-ark=arn} and, hence, 
with~$\ark(\chi)$. Applying the theorem in the introduction
completes the proof.
\end{proof}

\bibliography{references}
\bibliographystyle{plain}
\end{document}